\documentclass[a4paper]{amsart}

\usepackage{epsfig,amsmath}
\usepackage{xspace}
\usepackage[psamsfonts]{amssymb}
\usepackage[latin1]{inputenc}
\usepackage{float}\restylefloat{figure}

\usepackage{hyperref}

\newtheorem{proposition}{\bf Proposition}
\newtheorem{definition}[proposition]{\bf Definition}
\newtheorem{theorem}[proposition]{\bf Theorem}

\newtheorem{lemma}{\bf Lemma}
\newtheorem*{lemma*}{\bf Lemma}
\newtheorem*{sublemma*}{\bf Sublemma}
\newtheorem*{claim*}{\bf Claim}
\newtheorem*{complement*}{\bf Complement}
\newtheorem{corollary}[proposition]{\bf Corollary}
\newtheorem{conjecture}[proposition]{\bf Conjecture}

\renewcommand{\epsilon}{\varepsilon} 

\newcommand{\setof}[2]{\big\{{#1}\,\big|\,{#2}\big\}}

\def\tend{\longrightarrow}

\def\ds{\displaystyle}
\def\wt{\widetilde}

\def\on{\operatorname}

\def\ov{\overline}

\def\C{{\mathbb C}}
\def\D{{\mathbb D}}

\def\Q{{\mathbb Q}}

\def\R{{\mathbb R}}

\def\N{{\mathbb N}}

\def\cal{\mathcal}
\def\rad{\on{rad}}

\def\proof{\noindent{\bf Proof. }}

\def\ds{\displaystyle}

\newcommand{\avg}[1]{\underset{#1}{\on{avg}}}
\newcommand{\ninf}{_{n\to +\infty}}

\begin{document}

\title{A new proof of a conjecture of Yoccoz, Remarks, New results}
\author{Xavier Buff}
\author{Arnaud Chéritat}
\date{}

\begin{abstract}
We give a new proof of the following conjecture of Yoccoz:
\[(\exists C\in \R)~(\forall \theta\in \R)\quad \log R(P_\theta) \leq -Y(\theta) +C,\]
where $P_{\theta}(z)=e^{2i\pi \theta}z+z^2$, $R(P_\theta)$ is the conformal radius of the Siegel disk of $P_\theta$ (or $0$ if there is none) and $Y(\theta)$ is Yoccoz's Brjuno function.

In a former article we obtained a first proof based on the control of parabolic explosion.  Here, we present a more elementary proof based on Yoccoz's initial methods.

We then extend this result to some new families of polynomials such as $z^d+c$, $d>2$. We also show that the conjecture does not hold for $e^{2i\pi\theta} z + z^d$ with $d>2$.
\end{abstract}

\maketitle{}

\setcounter{tocdepth}{1}
\tableofcontents

\section{A new proof of a conjecture of Yoccoz}

\subsection{Introduction}

For a holomorphic germ $f$ of the form
\[f(z) = e^{2i\pi\theta} z + {\cal O}(z^2)\]
with $\theta\in\R\setminus\Q$, let
\[R(f) \in [0,+\infty]\]
be the radius of convergence of the linearizing formal power series, which we recall is the unique formal series of the form
\[h(Z)=Z+\sum_{n=2}^{+\infty} b_n Z^n\]
in $\C[[Z]]$ such that
\[f \circ h(Z) = h \circ R_\theta(Z)\]
with $R_\theta(Z)=e^{2i\pi\theta} Z$.
Let
\[P_\theta(z) = e^{2i\pi\theta} z + z^2.\]

\pagebreak

We do not recall the definition of Yoccoz's arithmetic Brjuno function
\[Y(\theta) = \sum_{n=0}^{+\infty} \theta_0 \cdots \theta_{n-1} \log \frac{1}{\theta_n}\]
where $\theta_0=\on{Frac}(\theta) = \theta - \lfloor \theta \rfloor$ and $\theta_{n+1} = \on{Frac}(1/\theta_n)$ when $\theta$ is irrational, and $Y(\theta)=+\infty$ if $\theta$ is rational\ldots\ Let's also not recall that a Brjuno number is an irrational real number $\theta$ satisfying Brjuno's condition:
\[\theta\in{\cal B} \iff Y(\theta)<+\infty.\]


\medskip

Yoccoz had proved in \cite{Y} that
\begin{theorem}[Yoccoz]\label{theo_phille}
  There exists a constant $C\in\R$ such that for all $\theta\in{\cal B}$, for all univalent function $f : \D \to \C$ fixing the origin with derivative $e^{2i\pi\theta}$,
  \[\log\, \on{inrad}(f) \geq -Y(\theta)-C.\]
  where $\on{inrad}(f)=\sup\setof{r\geq 0}{B(0,r)\text{ is contained in the Siegel disk of }f}$\footnote{the Siegel disk at the origin}.
\end{theorem}
\noindent Whence
\begin{corollary}[Yoccoz] $\forall \theta\in\cal B$,
  \[\log R(P_\theta)\geq - Y(\theta)-C - \log 2\]
\end{corollary}
\noindent the term $-\log 2$ coming from the fact $P_\theta$ is univalent only on the disk $B(0,1/2)$.

\medskip

He also had proved that
\begin{theorem}[Yoccoz]\label{theo_jasmin}
  There exists a constant $C\in\R$ such that for all $\theta\in\R$ there exists a univalent function $f : \D \to \C$ fixing the origin with derivative $e^{2i\pi\theta}$, such that
  \[\log R(f) \leq -Y(\theta)+C.\]
  This includes the case $\theta\notin\cal B$ (i.e.\ $Y(\theta)=+\infty$) if we interpret the above inequality as $R(f)=0$.
\end{theorem}
He then used a technique of Il'Yashenko and the polynomial-like map theory of Douady and Hubbard to transfer this result to $P_\theta$, but at the cost of a slight loss:
\begin{corollary}[Yoccoz]\label{cor_auxpieds}
  For all $\epsilon>0$, there exists $C_\epsilon\in\R$ (that a priori may tend to $+\infty$ as $\epsilon \tend 0$) such that for all $\theta\in\R$,
  \[\log R(P_\theta) \leq -(1-\epsilon)Y(\theta)+C_\epsilon.\]
  In particular, when $\theta$ is not Brjuno, $P_\theta$ is not linearizable.
\end{corollary}
So the natural conjecture was
\begin{conjecture}[Yoccoz]\label{conj_ecture}
  There exists a constant $C\in\R$ such that for all $\theta\in\R$,
  \[\log R(P_\theta) \leq -Y(\theta)+C.\]  
\end{conjecture}
\noindent (Since $Y$ is an unbounded function, this does not follow from corollary~\ref{cor_auxpieds}.)

\medskip

The second author found an independent\footnote{see footnote~\ref{foot} (when it appears)} proof of the optimality of Brjuno's condition in \cite{C}. It works directly in the family $P_{\theta}$. He looked at how parabolic points explode into cycles and how these cycles hinder each others. The control on parabolic explosion uses the combinatorics of quadratic polynomials, and the Yoccoz inequality on the limbs of the Mandelbrot set\footnote{\label{foot}Therefore this proof of optimality of the Brjuno condition is independent of Yoccoz's proof, but not of Yoccoz's overall work\ldots}.
The \emph{relative Schwarz lemma} of the first author then enabled us to have a good enough control on conformal radii to prove Yoccoz's conjecture in \cite{BC}:
\begin{theorem}[Buff, Chéritat] Conjecture~\ref{conj_ecture} holds.
\end{theorem}

The purpose of this article is to give a new proof of this conjecture, which turns out to be Yoccoz's proof of corollary~\ref{cor_auxpieds}, completed with a simple remark.

\begin{theorem}[Buff, Chéritat]\label{et_la_bete} Conjecture~\ref{conj_ecture} holds (with a different proof).
\end{theorem}

Following ideas of P\'erez-Marco~\cite{PM} and Geyer~\cite{G} we will also prove a new result: a similar statement holds for other families in which it was yet unknown\footnote{the technique of \cite{BC} turns out to be difficult to adapt in some of these new situations}. See section~\ref{sec_ond}. 

\subsection{The new proof}\label{sec_taire}\ \\

It follows closely Yoccoz's proof of corollary~\ref{cor_auxpieds}, but completes it.

His approach starts from any univalent function $f$ (for instance the one provided by theorem~\ref{theo_jasmin}).
Following Il'Yashenko, he considers the one-parameter family
\[f_a(z)=f(z)+az^2\]
for $a\in\C$. Then look at the formal linearizing power series of $f_a$:
\[h_a(Z)=Z + \sum_{n=2}^{+\infty} b_n(a) Z^n.\]

By Hadamard's theorem,
\[R(f_a) = 1/\limsup\ninf \sqrt[n]{|b_n(a)|}.\]

If we rescale by $1/a$: 
\[g_a(z) = a f_a(a^{-1}z) = af(a^{-1}z) + z^2.\]
This rescaled map tends to $P_\theta$ when $a\tend \infty$.
Moreover, we find in \cite{Y} page 60:
\begin{lemma}\label{lem_ing}
For all $f :\D \to \C$ univalent with $f(0)=0$ and $|f'(0)|=1$, for all $a\in\C$ with $|a|> 10$, the map $f_a$ has a quadratic-like restriction from $U$ to $V=B(0,13/36)$ (where $U=\setof{z\in B(0,1/3)}{f_a(z)\in V}$).
\end{lemma}

Now, we introduce in the argument the following two facts
(they both are valid in a much more general setting). Let $\avg{|a|=r}{m(a)}$ denote the average of the function $m(a)$ on the circle $|a|=r$ (with respect to the Lebesgue measure on the circle).

\begin{proposition}\label{prop_sac} If $\theta\in\cal B$, then
  $\log R(f) \geq \avg{|a|=10+1} \log R(f_a)$
\end{proposition}
\proof
As a holomorphic function of $a$, the function $b_n(a)$ satisfies
\[\log |b_n(0)| \leq \avg{|a|=10+1} \log |b_n(a)|.\]
We want to apply Fatou's lemma, but for this we need a uniform bound on the functions $b_n(a)$ for $|a|= 10+1$.

Proof n°1: there exists $r>0$ such that all functions $f_a$ are univalent on $B(0,r)$ for $|a|= 10+1$. By Yoccoz's theorem~\ref{theo_phille}, $\on{inrad}(f_a)$ is bounded from below for $|a|= 10+1$, say by $\eta>0$. Since $h_a$ is a univalent function on $D(0,\on{inrad} f_a)$, with derivative $1$ at $0$, this implies an upper bound of the form $\exists M\in\R$, $\forall n\geq 0$, $\forall a\in\C$ with $|a|= 10+1$, $\frac{1}{n} \log |b_n(a)| \leq M$.\footnote{We do not need the precise sharp value of the lower bound but just the existence of some for any fixed Brjuno number $\theta$. With this respect, Brjuno's original proof probably directly provides an upper bound on $\frac{1}{n} \log |b_n(a)|$, uniform on $|a|\leq 10+1$.}

Proof n°2: for $|a|=10+1$, we are in a domain ($|a|>10$) where the Siegel disk undergoes a holomorphic motion (see the proof of proposition~\ref{prop_ellant} for more details), thus this Siegel disk has its conformal radius bounded from below on the circle $|a|=10+1$, and thus $\frac{1}{n} \log |b_n(a)| \leq M$ holds for $|a|=10+1$.
\footnote{This second proof has the advantage of avoiding completely the arithmetic aspect, provided we replace the definition of $\cal B$ by $\theta\in{\cal B} \iff P_\theta$ is linearizable. The first proof works in a more general setting. See the remarks section.}

This uniform bound allows us to apply Fatou's lemma:
\[\limsup\ninf \frac{1}{n}\log |b_n(0)| \leq \avg{|a|=10+1} \limsup\ninf \frac{1}{n}\log |b_n(a)|.\]
i.e.\footnote{beware the minus sign}
\[\log R(f) \geq \avg{|a|=10+1} \log R(f_a).\]
\qed

The following is a corollary of lemma~\ref{lem_ing}:
\begin{proposition}\label{prop_ellant}
  The map $a\mapsto \log R(f_a)$ is harmonic for $|a|>10$.
\end{proposition}
\proof
  The quadratic-like restrictions all have a neutral fixed point. Therefore, this is the only non repelling cycle of the restriction. Therefore, the Julia set of the quadratic-like map undergoes a holomorphic motion as $a$ varies. The radius of convergence of $h_a$ coincides with the conformal radius of the Siegel disk $\Delta_a$ of the quadratic-like restriction. Now, when a Siegel disk has a boundary which undergoes a holomorphic motion, then its conformal radius has a logarithm $\log \rad \Delta_a$ that varies harmonically. [\small Let us give a proof of this fact which was communicated to us by Saeed Zakeri: first, note that the conformal radius varies continuously. Then, consider an extension\footnote{Slodkowsky's theorem provides one, but we can also use the Bers-Royden or the Sullivan-Thurston version since this argument is local in terms of the parameter.} of the holomorphic motion to a holomorphic motion of all the plane, but which does not necessarily commute with the dynamics. Let $a_0$ be any parameter and $z_n$ be any sequence in the Siegel disk of parameter $a_0$. For $a$ close to $a_0$ let $z_n(a)$ be the point that the motion transports $z_n$ to. Now look at $f_n(a)=h_a^{-1}(z_n(a))$. As the map $(a,z)\mapsto (h_a(z),z)$ is bi-analytic, $a\mapsto f_n(a)$ is analytic. Therefore, $\log |f_n(a)|$ is harmonic (it does not vanish). Now the map $a\mapsto \log \rad\Delta_a$ is the limit of these harmonic functions.\normalsize]
\qed

\medskip

Now,
\[\log R(g_a)=\log |a| +\log R(f_a)\]
is a harmonic function of $a \in \C\setminus \overline{B}(0,10)$, which tends to $R(P_\theta)$ when $a\tend \infty$. Therefore
\[\log R(P_\theta) = \avg{|a|=10+1} R(g_a) = \log 11 + \avg{|a|=10+1} \log R(f_a).\]
Together with proposition~\ref{prop_sac}
\[\log R(P_\theta) \leq \log 11 + \log R(f).\]
Then, using theorem~\ref{theo_jasmin} we get theorem~\ref{et_la_bete}.\footnote{whence an improvement of $4$}
\hfill Q.E.D.

\section{Remarks}

This section does not claim to bring new results. It is just a discussion of probably known and hopefully useful facts.

First, remember that for a germ $f=e^{2i\pi\theta} z + \ldots$, the radius of convergence $R(f)$ of its linearizing formal power series $h\in\C[[z]]$ is not necessarily equal to the conformal radius of the maximal linearization domain $\Delta(f)$ of $f$. An obvious possibility would be that $f$ has an extension to a bigger domain, which has a bigger maximal linearization domain. But this is not the only thing that can happen, since $h$ is not necessarily injective on its disk of convergence: in fact $h$ can be \emph{any convergent power series} of the form $z+{\cal O}(z^2)$. Indeed, for such an $h$, set $f$ as the conjugate near $0$ of $R_\theta$ by $h$\ldots

For instance, $h(z)= e^z-1=z+\ldots$ has infinite radius of convergence, and is not injective on $\C$.
The map $h(z)$ can also have critical points.

\bigskip\noindent\textsc{No arithmetics}\medskip

We can avoid completely the arithmetics discussion, by replacing Yoccoz's Brjuno function by
\[R(\theta) = \inf\setof{R(f)}{f\in{\cal S}_\theta}\]
where ${\cal S}_\theta$ is the set of univalent functions on $\D$ fixing $0$ with derivative $e^{2i\pi\theta}$.
Our enhancement of Yoccoz's method then takes the form:
\[(\exists C\in\R)\ (\forall \theta\in \R)\ \log R(P_\theta) \leq \log R(\theta) + C.\]

\bigskip\noindent\textsc{Superharmonicity of radii}\medskip

For all Brjuno number $\theta\in\cal B$, for all analytic family $f_\lambda(z) = e^{2i\pi\theta} z + \ldots$ of analytic germs, the function $\lambda\mapsto \log R(f_\lambda)$ is the opposite of a $\limsup$ of subharmonic functions:
\[\phi(\lambda)=\log R(f_\lambda) = -\limsup\ninf \frac{1}{n}\log|b_n(z)|.\]
By Yoccoz's theorem (\ref{theo_phille}) (or by Brjuno's computations), these functions are locally uniformly bounded above. Therefore, by Fatou's lemma, the function $\log R(f_\lambda)$ is (everywhere) above its average on circles.
But we can say more: by the Brelot-Cartan theorem (see \cite{Ra}), if we note $\phi^*(z)=\underset{z'\to z}{\liminf} \,\phi(z')$, then $\phi^*$ is superharmonic and $\phi=\phi^*$ except on a polar set.

We however cannot say that $\phi$ itself is superharmonic (iff $\phi=\phi^*$) because it is not necessarily lower semicontinuous, as the following counterexample shows.
Let $f_0=e^{2i\pi\theta} z +{\cal O}(z^2)$ be the restriction to $\D$ of a map $\wt{f}$ defined on an open set $\Omega$ containing $\ov\D$, and such that its Siegel disk $\Delta(\wt{f})$ in $\Omega$ goes over the edge of $\D$. For instance $\wt{f}=R_\theta(z)$ on $\Omega=\C$. Let $f_\lambda = f_0+\lambda z^2 g(z)$ for $\lambda\in\C$, where $g(z)$ is any analytic function on $\D$ that is singular on all of $\partial\D$. Then for $\lambda\neq 0$, the linearizing map $h_\lambda$ must map its disk of convergence in $\D$ (as in \cite{Y}). Also, and as a corollary, it is injective on its disk of convergence and maps it to the Siegel disk of $f_\lambda$. This implies that
\[ \underset{\lambda\underset{\neq}\to 0}\limsup\, R(f_\lambda) \leq \rad(U)< R(f_0)\]
where $U$ is the biggest $f_0$-invariant subdisk of $\Delta(\wt{f})$ that is contained in $\D$, and $\rad(U)$ is its conformal radius with respect to $0$.

Now, note that the lower semi-continuity holds if, instead of considering $R(f_\lambda)$ we consider $\rad(\Delta(f_\lambda\big|_\D))$ and if $f_\lambda \tend f_0$ for the compact open topology on $\D$. This is a corollary of the work of Risler \cite{Ris} (see also \cite{C} or \cite{BC36}).
Upper semi-continuity also holds, this time for an elementary reason: if $\rad \Delta(f_{\lambda_n})$ tends to some real $r$, the conformal maps $\phi_n$ make a normal family. They are also known to linearize $f_\lambda$. Therefore, any limit of the $\phi_n$ must linearize $f_\lambda$. Thus $\rad \Delta(f_\lambda) \geq r$. Analytic dependence on the parameter is not needed. Hence
\begin{proposition}
  For all $\theta\in\cal B$, the map $\ds\left(\begin{array}{rcl}H_\theta(\D) & \to & ]0,1] \\ f & \mapsto & \rad \Delta(f)\end{array}\right)$ is continuous, where $H_\theta(\D)$ is the set of analytic functions $f:\D\to\C$ fixing $0$ with derivative $e^{2i\pi\theta}$, and for the compact open topology\footnote{uniform convergence on compact subsets of $\D$} on $H_\theta(\D)$.
\end{proposition}

\begin{proposition}
  If $U$ is a one complex dimensional parameter space and $(\lambda,z)\in U\times \D\mapsto f_\lambda(z) = e^{2i\pi\theta} z + {\cal O}(z^2)$ is analytic, then the map
  \[\lambda \mapsto \log \rad \Delta(f_\lambda)\]
  is continuous and superharmonic.
\end{proposition}
\begin{proof}
We already mentioned the continuity.
\\
Now, here is a trick\footnote{It would be nice to have a more satisfactory (no power series) proof. Also, it could be true that superharmonicity still holds if the domain of definition of $f$ undergoes a holomorphic motion.} that yields superharmonicity with little effort:
consider a function $g$ as in the discussion above, i.e.\ holomorphic on $\D$ and with singularities at all points of $\partial\D$. Consider the sequence of families $(\lambda,z)\mapsto \tau_n^{-1} f_\lambda(\tau_n z) +\frac{1}{n}z^2g(z)$ with $\tau_n=1-\frac{1}{n}$. They all satisfy $\log R=\log \rad \Delta$ (as in~\cite{Y}), whence all these are (continuous) superharmonic functions of $\lambda$. By the previous proposition, these functions tend (locally uniformly) to $\log \rad \Delta(f_\lambda)$.
\end{proof}

\bigskip\noindent\textsc{Harmonicity}\medskip

\begin{proposition}
  If $f_\lambda : U_\lambda \to \C$ is an analytic family of maps of the form $f(z) = e^{2i\pi\theta} z + {\cal O}(z^2)$ and if the boundary of $\Delta(f_\lambda)$ undergoes a holomorphic motion (we do not require $\partial \Delta(f_\lambda) \subset U_\lambda$), then the function $\lambda \mapsto \log \rad \Delta(f_\lambda)$ is harmonic.
\end{proposition}
\begin{proof}
  Same as in the second part of the proof of proposition~\ref{prop_ellant} (courtesy of S.\ Zakeri).
\end{proof}

\bigskip\noindent\textsc{Superharmonicity under holomorphic motions}\medskip

Let us also mention the following
\begin{proposition}
  If a simply connected open subset of $\C$ undergoes a holomorphic motion (of its boundary) then the logarithm of its conformal radius with respect to a holomorphically varying point in it is a (continuous) superharmonic function of the parameter.
\end{proposition}
\begin{proof}
  Let $U_\lambda$ be this set, $c_\lambda$ be the varying point and $r(\lambda)$ be the conformal radius of $U_\lambda$ w.r.t.\ $c_\lambda$. Let $V_\lambda$ be the image of $U_\lambda$ by the inversion $z\mapsto 1/(z-c_\lambda)$. The set $V_\lambda$ is bounded and undergoes a holomorphic motion of its boundary. The conformal radius of $U_\lambda$ is the inverse of the capacity radius of $V_\lambda$, which is itself expressable by an energy minimization as follows:\footnote{As a variant of this, one could instead use the transfinite diameter.}
  \[\log r(\lambda) = -\log \text{capacity radius} = \inf_\mu E(\mu)\]
  where $\mu$ varies in the set of probability measures on $\partial V_\lambda$ and $E(\mu)$ (the energy) is defined by
  \[E(\mu) = \int_{\partial V_\lambda\times\partial V_\lambda} -\log \big|u-v\big| d\mu(u)d\mu(v)\]
  where the integrand is understood to be $+\infty$ when $u=v$.
  Choose a basepoint $\lambda_0$ and let $\xi_\lambda(z) : \partial V_{\lambda_0} \to \partial V_\lambda$ be the holomorphic motion.
  Then, for all probability measure $\mu$ on $\partial V_{\lambda_0}$,
  \[E\big((\xi_\lambda)_*\mu\big) = \int_{\partial V_{\lambda_0}\times\partial V_{\lambda_0}} -\log \big|\xi_\lambda(u)-\xi_\lambda(v)\big| d\mu(u)d\mu(v).\]
  This is a harmonic function of $\lambda$. Now, the infimum in the energetic definition of $r(\lambda)$ yields a superharmonic function.
\end{proof}

\bigskip\noindent\textsc{How strange}\medskip

This is kind of surprising: let $A$ denote the fact that a simply connected domain undergoes a holomorphic motion (of its boundary), and $B$ denote the fact that this domain is a Siegel disk of an analytically varying family of analytic maps (with fixed rotation number) in $\cal \D$. Then 
\begin{center}\parbox{6cm}{
$A \implies$ $\log \rad$ is superharmonic,
\par\smallskip
$B \implies$ $\log \rad$ is superharmonic,
\par\smallskip
$(A \text{ and } B) \implies$ $\log \rad$ is harmonic\ldots
}\end{center}
Is it fair that when a number has two reasons to be negative, then it is null?

\bigskip\noindent\textsc{Other radii of interest}\medskip

We have
\[R(f) = \text{ the radius of convergence of }h\]
and
\[\rad \Delta(f) = \text{ the biggest radius }\leq R\text{ below which }h\text{ maps in }\Delta(f).\]
Here are a few other ``natural'' radii that one could study
\begin{eqnarray*}
A &=& \text{ the biggest radius }\leq R\text{ on which }h\text{ is injective,}\\
B &=& \text{ the biggest radius }\leq R\text{ on which }h\text{ has no critical point,}\\
C &=& \text{ the biggest radius }\geq R\text{ on which }h\text{ has a meromorphic extension }\wt{h},\\
D &=& \text{ the biggest radius }\leq C\text{ on which }\wt{h}\text{ is injective,}\\
E &=& \text{ the biggest radius }\leq C\text{ on which }\wt{h}\text{ has no critical point.}
\end{eqnarray*}

\section{New results}\label{sec_ond}

Let us define a class of well behaved polynomials that was studied by Lukas Geyer in~\cite{G}. An orbit tail is an equivalence class in the set of forward critical orbits\footnote{or, equivalently, in the set of critical \emph{points}}, with the relation $z\equiv z' \iff \exists m,n \in\N$ such that $P^n(z)=P^m(z')$. We say it is infinite if a point of the class (and therefore every point in the class) has infinite forward orbit.
\begin{definition}
 We will say that a polynomial has ``the above property'' if the number of infinite critical orbit tails is equal to the number\footnote{it is at least this number} of indifferent cycles.
\end{definition}

\begin{remark} If $f$ satisfies "the above property", all its iterates do not necessarily.
\end{remark}

Lukas Geyer proved optimality of Brjuno's condition for these polynomials (and even for a bigger class\footnote{for the class of \emph{saturated} polynomials, i.e.\ polynomials such that the number of infinite critical orbit tails \emph{in the Julia set} is equal to the number of indifferent cycles}), by using the same method as Yoccoz. It is therefore natural that our new observation adapts in this setting:

\begin{definition}\label{ilenmer}
 Critical orbits are the sets $\setof{P^k(c)}{k\geq 0}$ where $c$ is a critical point of $P$.
 A point $z$ in a critical orbit is said to be free\footnote{this is not a standard terminology} if for all critical point $c'$, $\forall k\in\N$, $\forall l\in\N^*$,
 $P^k(c')=P^l(z) \implies$ $k\geq l$ and $P^{k-l}(c')=z$.
 We will say that a polynomial has ``the above property with bound $N$'' if moreover the cardinal of the union of all indifferent cycles and the set $Z$ is $\leq N$, where $Z$ is the set of non-free points of critical orbits.
\end{definition}

\begin{theorem}\label{theo_dore}
  Let $N\in\N$ and $\cal C$ be a compact set of degree $d$ polynomials $f$ fixing $0$ with indifferent multiplier $e^{2i\pi\theta(f)}$, with ``the above property with bound $N$'', and moreover such that the indifferent cycles stay bounded away from $0$ when $f$ varies in $\cal C$.\footnote{It turns out that this last hypothesis is automatic if we require the other ones. The proof uses the notion of virtually repelling parabolic point. However, in our applications, it does not cost much to check it directly. This relieves the proof of theorem~\ref{theo_dore}.}
  Then $\exists C\in\R$ such that $\forall f\in \cal C$,
  \[\log R(f) \leq -Y(\theta(f))+C.\]
  Hence, by Yoccoz's theorem~\ref{theo_phille}, 
  \[\Big|\log R(f) + Y(\theta(f))\Big| \text{ is bounded over }{\cal C}.\]
\end{theorem}
\proof
Let \[G_f(z)=\prod_i(z-w_i)^{n_i}\prod_j(z-u_j)^2\]
where $\{w_i\}=Z\setminus\{0\}$ ($Z$ was defined in definition~\ref{ilenmer}), $n_i =$ the local degree of $f$ at $w_i$, and $u_j$ are the indifferent periodic points of $f$ (including $0$).
Let \[g_{f,a}=f + a G_f.\]
By compactness of $\cal C$, by the bound $N$, and by the definition of $G_f$, we see that (a)~$G_f$ is a bounded family over $\cal C$, and moreover by the last hypothesis that (b)~$G_f''(0)$ is bounded away from $0$.
Therefore, by (a) there exists $r>0$ and $R>0$, independent of $f$, such that $|a|<r \implies g_{f,a}(z)$ is a polynomial-like map of degree $d$ from the-component-of-$g_{f,a}^{-1}(B(0,R))$-contained-in-$B(0,R)$ to $B(0,R)$. For any fixed $f\in\cal C$, as $a$ varies in $B(0,r)$, this polynomial-like map cannot undergo a parabolic bifurcation. Therefore its Julia set undergoes a holomorphic motion. Therefore the same analysis as in section~\ref{sec_taire} holds (with parameter $1/a$ instead of $a$, but this is completely equivalent) and we can thus write:
\[\forall f\in{\cal C},\quad\log R(f) = \avg{|a|=r/2} \log R(g_{f,a}).\]
Let $\wt{g}_{f,a} (z)= a g_{f,a} (a^{-1} z)$. Then, as $a\tend \infty$, $\wt{g}_{f,a}$ tends (pointwise) to the degree~$2$ polynomial 
\[P(z)= e^{2i\pi\theta(f)} z + \frac{G_f''(0)}{2} z^2\]
The same analysis as in section~\ref{sec_taire} also holds and yields
\[\log R(P) \geq \avg{|a|=r/2} R(\wt{g}_{f,a}).\]
Now
\[\log R(P) = \log R(e^{2i\pi\theta(f)} z + z^2) - \log \frac{|G_f''(0)|}{2}\]
and
\[\log R(\wt{g}_{f,a}) = \log R(g_{f,a}) + \log|a|\]
and we proved that
\[\log R(e^{2i\pi\theta(f)} z + z^2) \leq -Y(\theta(f)) + C.\]
Putting it altogether, we get
\[\log R(f) \leq -\log{r/2} - \log \frac{|G_f''(0)|}{2} - Y(\theta(f)) + C.\]
By (b), we get the upper bound of the theorem.
To apply then Yoccoz's theorem, there remains to remark that the maps $f\in\cal C$ are all univalent on a common disk $B(0,r')$, since they have bounded degree and critical points are necessarily bounded away from $0$.
\qed

\begin{remark} We did not try to get the most general result possible. For instance, it is possible that the hypothesis that $0$ has period $1$ is not required.\footnote{But in this case, there is one more condition: that the indifferent cycle $0$ belongs to does not collapse on itself. This condition is not implied by the others.}
\end{remark}


\begin{corollary}
  This holds for the boundary of the central hyperbolic component of the family $z^d+c$, ($d\geq 2$).
\end{corollary}

\begin{corollary}
  This holds also for the family $e^{2i\pi\theta}z(1-z)^{d-1}$.
\end{corollary}
\noindent Indeed, the critical points are $z=1/d$ and $z=1$ (with multiplicity $d-2$). The second critical point is mapped in one step on $z=0$. Thus we may apply theorem~\ref{theo_dore} with $N=2$.

\begin{corollary}
  For the family
  \[f_\theta(z)=e^{2i\pi\theta}(z + z^d),\]
  the following holds: $\exists C>0$ such that $\forall \theta\in\R$,
  \[-\frac{Y\big((d-1)\theta\big)}{d-1}-C \leq \log R(f_\theta) \leq -\frac{Y\big((d-1)\theta\big)}{d-1}+C.\]
  As a consequence, the function
  \[\theta\in{\cal B} \mapsto \log R(f_\theta)+Y(\theta)\]
  is unbounded (on any interval). On the other hand, 
  \[-Y(\theta)-C' \leq \log R(f_\theta) \leq -\frac{Y(\theta)}{d-1}+C'.\]
\end{corollary}
\proof
  The family $f_\theta$ is semi-conjugated to the previous family: more precisely let
  $\phi(z)=-z^{d-1}$, and $g_\theta(z)=e^{2i\pi\theta}z(1-z)^{d-1}$. Then\footnote{Note how the rotation number changed.}
  \[\phi\circ f_\theta = g_{(d-1) \theta} \circ \phi.\]
  The first claim follows at once.
  The second claim (unboundedness) follows directly from the fact that $Y$ is an unbounded function. For the last claim, the left inequality follows from Yoccoz's theorem~\ref{theo_phille},\footnote{Whence a non-arithmetic proof of the following fact: $\forall m$, $\exists C>0$ $\forall \theta$, $Y(m \theta) \leq m Y(\theta) + C$. This fact can also be proved directly from the definition of $Y$.} the right inequality follows from the following arithmetic lemma.\footnote{we need a little bit weaker: $\forall m\in\N^*$, $\exists C_m\in\R$,  $\forall \theta\in\R$, $Y(\theta) \leq Y(m\theta) + C_m$} \qed
\begin{lemma}
  $\exists C>0$, $\forall m\in\N^*$, $\forall \theta\in\R$, \[Y(\theta) \leq Y(m\theta) + C\log m.\]
\end{lemma}
\proof
  We will use the Brjuno sum:
  \[B(\theta) = \sum_{n\in\N} \frac{\log q_{n+1}}{q_n}\]
  where $p_n/q_n$ are the approximants of $\theta$.
  We have the following arithmetical property (c.f.\ \cite{Y}, page~14):
  \[\big|B(\theta)-Y(\theta)\big|\text{ is bounded}.\]
  We recall that (a) if $p_n/q_n$ are the approximants of $\alpha$ then $1/2q_n q_{n+1}<|\alpha-p_n/q_n|<1/q_n q_{n+1}$, and also that (b) $q_n \geq F_n$ where $F_n$ is the $n$-th Fibonacci number. Last, (c) if $|\alpha-p/q|<1/2q^2$, then $p/q$ is an approximant of $\alpha$.\\
Now, for every approximant $p_n/q_n$ of $\theta$, note $m p_n/q_n = p'/q'$ with $q' = q_n / (m \wedge q_n)$.\\
Either $p'/q'$ is itself an approximant of $m\theta$ in which case if we note $p''/q''$ the next approximant of $m\theta$, then $1/2q'q''<|m\theta-p'/q'|=m|\theta-p_n/q_n| < m/q_n q_{n+1}$ whence $q''> q_{n+1}q_n/2mq' > q_{n+1}/2m$ and thus $\frac{\log q''}{q'} \geq \frac{\log q_{n+1}}{q'} - \frac{\log 2m}{q'} \geq \frac{\log q_{n+1}}{q_n} - \frac{\log 2m}{q'}$. 
Or $m p_n/q_n = p'/q'$ is not an approximant of $m\theta$, which means that $|m\theta - p'/q'| \geq \frac{1}{2q'^2}$, and thus
$\frac{1}{q_n q_{n+1}}\geq |\theta - p_n/q_n|\geq \frac{1}{2mq'^2}$ whence $q_{n+1}\leq \frac{2mq'^2}{q_n} \leq 2mq_n$ and thus $\frac{\log q_{n+1}}{q_n} \leq \frac{\log q_n}{q_n} + \frac{\log 2 m}{q_n}$. Finally,
\begin{eqnarray*}
B(\theta) & = & \sum_{\text{case }1} \frac{\log q_{n+1}}{q_n} + \sum_{\text{case }2} \frac{\log q_{n+1}}{q_n} \\
  & \leq & \sum \frac{\log q''}{q'} + \log(2m) \sum \frac{1}{q'} + {\sum}' \frac{\log F_n}{F_n} + \log (2m) \sum \frac{1}{F_n}
\end{eqnarray*}
the prime in the sum means the summand needs to be replaced by the smallest non increasing sequence greater or equal to the sequence $\log F_n/F_n$. 
For different values of $n$, the approximants $p'/q'$ of $m\theta$ are different since $p'/q' = m p_n/q_n$ and thus
\begin{eqnarray*}
B(\theta) & \leq & B(m\theta) + \log(2m) \sum \frac{1}{F_n} + {\sum}' \frac{\log F_n}{F_n} + \log (2m) \sum \frac{1}{F_n}
\end{eqnarray*}
Since $F_n$ is exponentially increasing, the sums (independent of $\theta$) they are involved in are finite.
\qed
\medskip

This suggests the following
\begin{conjecture}
  There exists a $C=C(d)\in\R$ such that for all polynomial $f$ of degree $d$ with an indifferent fixed point at the origin,
  \[\log R(f) \leq -\frac{Y(\theta)}{d-1}+\log \min|c_i|+C\]
  where the $c_i$ are the critical points of $f$ and $\theta$ is the rotation number at the origin.
\end{conjecture}
\noindent with possible refinements according to how many recurrent critical points are associated to the indifferent fixed point.

\end{document}